\documentclass[12pt]{article}
\usepackage[cp1251]{inputenc}
\usepackage[english]{babel}
\usepackage{amssymb,amsmath}
\newtheorem{Def}{Definition}[section]
\newtheorem{Ex}{Example}[section]
\newtheorem{Th}{Theorem}[section]

\newtheorem{Lem}{Lemma}[section]
\newtheorem{Prop}{Proposition}[section]
\newtheorem{Cor}{Corollary}[section]
\newtheorem{Rem}{Remark}[section]
\newenvironment{Proof}
{\par\noindent{\bf Proof.}}
{\hfill$\scriptstyle\blacksquare$}

\title{Order-Preserving Variants of The Basic Principles of Functional Analysis}
\author{A. A. Zaitov\footnote {Tashkent institute of architecture and civil engineering}}

\begin{document}

\maketitle
\thispagestyle{empty}

\begin{center}

\end{center}
\begin{abstract}
We established order-preserving versions of the basic principles  of functional analysis such as Hahn-Banach, Banach-Steinhaus, open mapping and Banach-Alaoglu theorems.\\


\end{abstract}

\tableofcontents

\section{Introduction}

Recently researches in the field of idempotent mathematics and also Choquet integrals intensively develop. Since its introduction in 1974 by Sugeno, the concept of fuzzy measure has been often used in multicriteria decision making. Later in \cite{grabroub2000appChoquetint}, the authors explained the methodology of using the Choquet integral in multicriteria decision making. The notion of idempotent measure (Maslov integral) finds important applications in different part of mathematics, fuzzy topology, mathematical physics and economics (see the article \cite{radul2018idemmeasurAR} and the bibliography therein). As well known idempotent measures and Choquet integrals are weakly additive, order-preserving functionals. But for this functionals there not establish yet the basic principles (analogous principles of Functional Analysis). In the present paper we will establish order-preserving versions of the the basic principles of Functional Analysis such as the Hahn-Banach, Banach-Steinhaus, open mapping and Banach-Alaoglu theorems.

Remind that \textit{partially ordered vector space} is a pair $(L, \leqslant)$ where $L$ is a vector space over the field $\mathbb{R}$ of real numbers, $\leqslant$ is an order satisfying the following conditions:

$1)$ if $x\leqslant y$, then $x+u\leqslant y+u$ for all $x,\ y,\ u\in L$;

$2)$ if $x\leqslant y$, then $\lambda x\leqslant \lambda y$ for all $x,\ y\in L$ and $\lambda\in \mathbb{R}_+$.

If the conditions $1)$ and $2)$ hold then they say that $\leqslant$ is linear order. A formation of a vector space $L$ with linear order $\leqslant$ over $\mathbb{R}$ is equivalent to indicate a set $L_+\subset L$ called a \textit{positive cone} in $L$ and owning the properties:
\[
L_++L_+\subset L_+;\qquad \lambda L_+\subset L_+,\quad \lambda\in \mathbb{R}_+; \qquad L_+\cap (-L_+)=\{0\}.
\]

In this case the order $\leqslant$ and the positive cone $K$ are connected by a relation
\[
x\leqslant y\Leftrightarrow y-x\in L_+, \qquad x,\ y\in L.
\]
Elements of $L_+$ is called \textit{positive vectors} of $L$.

Let $(L, L_+)$ be a partially ordered vector space. We say \cite{vern2009vectspordunit} that
$L_+$ \textit{is full} (or that $L_+$ \textit{is a full cone}) if $L=L_+-L_+$.

Let $x\in L_+$. The point $x$ is said to be \textit{an inner point} of the cone $L_+$ if for any segment $[x_1, x_2]$ containing $x$ as an inner point, the segment $[x_1, x_2]\cap L_+$ also contains it as an inner point. The set of all inner points of the cone $L_+$ is called \textit{an interior} of this cone, and it denotes as $Int L_+$.

Fix an inner point $x_0\in L_+$. For a $\delta> 0$ we determine a $\delta$-\textit{neighbourhood} (with respect to the cone $L_+$ and the point $x_0$) of zero $0\in L$ as following:

\numberwithin{equation}{section}
\begin{equation}\label{0neigh}
\langle0; \delta \rangle=\{x\in L\colon \ (\delta x_0\pm x)\in Int L_+\}.
\end{equation}
It is easy to see that a family of the sets of the view (\ref{0neigh}) forms a base of neighbourhoods of zero. A neighbourhood of an arbitrary point $z\in L$ can be defined by the shifts of the neighbourhoods of zero:
\begin{multline}\label{neigh}
\langle z; \delta \rangle =\langle0; \delta \rangle + z= \{x+z \in L\colon \ x\in \langle0; \delta \rangle\}=\\
=\{x+z \in L\colon \ (\delta x_0\pm x)\in Int L_+\} =\\
= \{y\in L\colon \ (\delta x_0\pm (y-z))\in Int L_+\}.
\end{multline}

\begin{Prop}\label{basetop}
A collection
\[
\left\{\langle z; \delta \rangle\colon \ z\in L,\ \delta>0 \right\}
\]
forms a base of a Hausdorff topology on $L$. Further, $L$ equipped with this topology becomes a topological vector space.
\end{Prop}

\begin{Proof} The proof consists of direct checking.

\end{Proof}

An element $1\in L$ of a partially ordered vector space $L$ is called \textit{(strongly) order unit} if $L=\bigcup\limits_{n=1}^{\infty}[-n 1,\ n 1]$. This is equivalent to what for every $x\in L$ there exists  $\lambda>0$ such that $-\lambda 1\leqslant x\leqslant \lambda 1$.
Let $x\in L$. A partially ordered vector space $L$ is called \textit{Archimedean} if the inequality $nx\leqslant 1$ executed for all $n=1,\ 2, \dots$, implies $x\leqslant 0$. In this case on $L$ one can define a norm by the equality
 \begin{equation} \label{norm}
\|x\|=\inf\{\lambda>0\colon \ -\lambda 1\leqslant x \leqslant \lambda 1 \}.
\end{equation}

The obtained norm is said to be \textit{an order norm}. A partially ordered vector space $L$ is called \textit{a vector space with an order unit} if $L$ has an order unit and $L$ is an Archimedean space. A topology on $L$ generated by the norm (\ref{norm}) is called \textit{order (vector)} topology. For a subset $X\subset L$ by $Int X$ we denote the interior of $X$ according to the order topology on $L$. We accept the following agreement
\[
x<y\Leftrightarrow y-x\in Int L_+.
\]

A set $U(0_E, \varepsilon) = \{x\in E: -\varepsilon 1_E<x<\varepsilon 1_E\}$ is an open neighbourhood of zero $0_E$ concerning to the order topology. As vector topology is invariant according to the shifts then for every point $x\in E$ a set $U(x, \varepsilon) = \{y\in E: -\varepsilon 1_E<y-x<\varepsilon 1_E\}$ is an open neighbourhood of $x$ with respect to the order topology.

\begin{Prop}
The order topology and topology introduced by Proposition \ref{basetop} on a vector space with an order unit coinside.
\end{Prop}

The Proof is trivial.

\section{Extensions of Order-Preserving Functionals}

In this section we will prove the order-preserving functional's variant of the Hahn-Banach theorem, one of the basic principles of functional analysis.

Let $L$ be a partially ordered vector space over the field $\mathbb{R}$ of real numbers, and $L_+$ be a full cone in it. Let $x_1$, $x_2\in L$ be arbitrary various points. The set $[x_1, x_2]=\{\alpha x_1+ (1-\alpha)x_2\colon \alpha\in [0, 1]\}$ is called \textit{a segment} connecting points $x_1$ and $x_2$. A point $x\in [x_1, x_2]$ is \textit{an inner point} of the segment $[x_1, x_2]$ if $x_1\neq x\neq x_2$.

\begin{Def}\label{Asubspace}
A subset $B$ of a partially ordered vector space $L$ is said to be an $A$-subspace concerning a point $x_0\in L$ if $0 \in B$, and $x\in B$ implies $(x+\lambda x_0)\in B$ for each $\lambda\in \mathbb{R}$.
\end{Def}

The following assertion is evident.

\begin{Lem}
A subspace $B$ of the partially ordered vector space $L$ is an $A$-subspace according to $x_0$ iff it contains $x_0$.
\end{Lem}

Note that the space $L$ and its subspace $\{\lambda x_0\colon \ \lambda\in \mathbb{R}\}$ are \textit{trivial} $A$-subspaces. As distinct from linear case the set $\{0\}$ is not $A$-subspace.

It is easy to see that an intersection of any collection of $A$-subspaces is a $A$-subspace. In particular, an intersection of all $A$-subspaces containing a given set $X$ is the minimal $A$-subspace, containing $X$; this $A$-subspace we call as \textit{a weakly additive span} of $X$, and designate through $A(X)$. The following statement describes a structure of the weakly additive span of a given set.

\begin{Prop}\label{spanstructure}
A weakly additive span $A(X)$ of a subset $X$ of a partially ordered linear  space $L$ consists of a (set-theoretic) union of $\{\lambda x_0\colon \ \lambda\in \mathbb{R}\}$ and the collection of all sums of the look $x+\lambda x_0$, $x\in X$, $\lambda\in \mathbb{R}$, i.~e.
\[
A(X)=\{\lambda x_0\colon \ \lambda\in \mathbb{R}\}\cup \bigcup\limits_{\substack {x\in X,\\ \lambda\in \mathbb{R}}}\{x+\lambda x_0\} = \bigcup\limits_{\substack {x\in X\cup\{0\},\\ \lambda\in \mathbb{R}}}\{x+\lambda x_0\},
\]
in particular, if $x_0\in X$ then
\[
A(X)=\bigcup\limits_{\substack {x\in X,\\ \lambda\in \mathbb{R}}}\{x+\lambda x_0\}.
\]
\end{Prop}

The proof is obvious.

Let's denote
\[
\Lambda=\{\lambda x_0;\ \lambda\in \mathbb{R}\}.
\]
Then we have
\[
A(X)=\bigcup\limits_{x\in X\cup\{0\}}(x+\Lambda).
\]

The last equality explains the name `$A$-subspace'.  Every $A$-subspace $A(X)$ consists of the union of one-dimensional subspace $\Lambda\subset L$ and affine subsets $x+\Lambda\subset L$, $x\in X$.

\begin{Def}\label{weakordnorm}
A functional $f\colon L\to \mathbb{R}$ is called:

$1)$ weakly additive (according to the point $x_0$) if
\[
f(x+\lambda x_0)=f(x)+\lambda f(x_0), \qquad x\in L,\quad \lambda\in \mathbb{R};
\]

$2)$ order-preserving (concerning to the cone $K$) if for every pair $x$, $y\in L$ belonging $y-x\in K$ implies the inequality
\[
f(x)\leqslant f(y);
\]

$3)$ normed (with respect to the point $x_0$) if $f(x_0)=1$.
\end{Def}

From the definition immediately follows that weakly additive functional is linear on the one-dimensional subspace $\{\lambda x_0;\ \lambda\in \mathbb{R}\}$ of $L$. From here we have $f(0)=0$.

Let $(L, L_+)$ be an partially ordered real vector space. A functional
$f\colon L\to \mathbb{R}$ is called \textit{positive} if $f(L_+) \subseteq [0, +\infty)$. Each weakly additive, order-preserving functional is positive. Really, let $x\in K$. Then $x-0\in K$. Since $f$ is order-preserving functional, then $f(x)\geqslant f(0)$. Consequently, $f(x)\geqslant 0$. There exists a functional which is weakly additive, positive but does not order-preserving.

\begin{Ex}\label{positbutnotord}
{\rm Let $L=\mathbb{R}^2=\{(x_1, x_2)\colon \ x_i\in \mathbb{R},\ i=1,\ 2\}$ be partially ordered vector space with respect to the usual linear operations `$\cdot$' -- the multiplication by real numbers, `$+$' -- the sum of elements of $L$, and to the pointwise order $\leqslant$ on $L$, which defines as $(x_1, x_2)\leqslant(y_1, y_2)\Leftrightarrow x_1\leqslant y_1 \mbox{ and } x_2\leqslant y_2$. The set $L_+=\{(x_1, x_2)\in \mathbb{R}^2\colon \ x_i\geqslant 0, i=1,\ 2\}$ is a positive cone in $L$. Define a functional $f\colon \mathbb{R}^2\to \mathbb{R}$ by the rule
\[
f(x_1, x_2)=\frac{1}{2}\left(x_1+x_2+\sqrt{|x_2-x_1|}\right), \qquad (x_1, x_2)\in \mathbb{R}^2.
\]
It is clear that $f$ is a weakly additive and positive functional. But we have $f\left(\frac{1}{2}, \frac{1}{2}\right)<f\left(\frac{1}{4}, \frac{1}{2}\right)$ though $\left(\frac{1}{2}, \frac{1}{2}\right)-\left(\frac{1}{4}, \frac{1}{2}\right)\in K$.}
\end{Ex}

\begin{Prop}\label{continuous}
If an order-preserving, weakly additive functional $f\colon L\to \mathbb{R}$ is continuous at zero $0$ then it is continuous on all $L$.
\end{Prop}

\begin{Proof}
Let for every $\varepsilon>0$ there exist $\delta>0$ such that $|f(x)|<\varepsilon$ for all $x\in \langle0; \delta\rangle\subset L$. Let $y\in L$ be an arbitrary nonzero element. Consider a neighbourhood
\[
\left\langle y; \frac{\delta}{2}\right\rangle= \left\{z\in L\colon \ \left(\frac{\delta}{2}x_0\pm (z-y)\right)\in Int L_+\right\}.
\]

For every $z\in \left\langle y; \frac{\delta}{2}\right\rangle$ we have:

1) $f(\frac{\delta}{2} x_0)+f(y)\geqslant f(z)$ since $\left(\frac{\delta}{2}x_0+y\right)-z\in Int L_+$;

2) $f(\frac{\delta}{2} x_0)+f(z)\geqslant f(y)$ since $\left(\frac{\delta}{2}x_0+z\right)-y\in Int L_+$.

From here follows that
\begin{equation}\label{absvalless}
|f(z)-f(y)|\leqslant f\left(\frac{\delta}{2} x_0\right).
\end{equation}
Now we take $x\in L$ such that $x-\frac{\delta}{2}x_0\in L_+$ and $\frac{3\delta}{4}x_0-x\in L_+$. Then $f(\frac{\delta}{2}x_0)\leqslant f(x)$ and $f(x)\leqslant f(\frac{3\delta}{4}x_0)$. On the other hand $f(\frac{3\delta}{4}x_0)<\varepsilon$ so far as $\frac{3\delta}{4}x_0\in \left\langle 0; \delta\right\rangle$. Consequently $|f(z)-f(y)|<\varepsilon$ for each $z\in \left\langle y; \frac{\delta}{2}\right\rangle$. So $f$ is continuous at $y\in L$. Thus $f$ is continuous on all $L$ owing to arbitrariness of $y\in L$.

\end{Proof}

A weakly additive, order-preserving functional $f\colon L\to \mathbb{R}$ is called \textit{bounded} if $\sup\{|f(x)|\colon \ x\in \langle 0; 1\rangle\}<\infty$.

\begin{Prop}\label{bounded}
A weakly additive, order-preserving functional is  bounded if and only if it is continuous.
\end{Prop}

\begin{Proof}
Let $f\colon L\to \mathbb{R}$ be weakly additive, order-preserving bounded functional. Let $f(x_0)=a<\infty$ and $\delta x_0 \pm(z-y)\in Int L_+$. Then similarly to (\ref{absvalless}) one can show that $|f(z)-f(y)|< \delta f(x_0)=\delta a$, and consequently $f$ is continuous.

Conversely, let a weakly additive, order-preserving functional $f\colon L\to \mathbb{R}$ be continuous. Then there exists $\delta> 0$ such that $|f(x)|<1$ at all $x\in \langle 0; \delta\rangle$. In particular, $\left|\frac{\delta}{2}f(x_0)\right|<1$ so far as $\frac{\delta}{2}x_0\in \langle 0; \delta\rangle$. Hence, $|f(x_0)|<\frac{2}{\delta}<\infty$, i. e. $\sup\{|f(x)|\colon \ x\in \langle0; 1\rangle\}<\frac{2}{\delta}<\infty$.

\end{Proof}

\begin{Cor}
A weakly additive, order-preserving, normed functional is continuous (or, the same, bounded).
\end{Cor}

The following statement is an analog of Hahn-Banach theorem for weakly additive, order-preserving functionals.

\begin{Th}\label{HahnBanach}
Let $B$ be an $A$-subspace of the space $L$. Then for every weakly additive, order-preserving functional $f\colon B\to \mathbb{R}$ there exists a weakly additive, order-preserving functional $f_0\colon L\to \mathbb{R}$ such that $f_0|_{B}=f$.
\end{Th}

\begin{Proof}
Let $y\in L\setminus B$. Put $B^{\prime}=B\cup \{y+\lambda x_0\colon \ \lambda\in \mathbb{R}\}$. Obviously that $B^{\prime}$ is an $A$-subspace of $L$. Put
\[
B^+=\{z\in B\colon \ z-y\in B\}\qquad \mbox{and} \qquad B^-=\{z\in B\colon \ y-z\in B\}.
\]
The obtained sets $B^+$ and $B^-$ are not empty. Indeed, take $\lambda>0$ such that $y\in \langle0; \lambda\rangle$. Then evidently that $2\lambda x_0\in B^+$ and $-2\lambda x_0\in B^-$.

Put
\[
p^+=\inf\{f(z)\colon \ z\in B^+\}, \qquad p^-=\sup\{f(z)\colon \ z\in B^-\}.
\]
We have $p^-\leqslant p^+$. Indeed, $y-z_1\in L_+$ if provided $z_1\in B^-$, and $z_2-y\in L_+$ if provided $z_2\in B^+$. From here we get $z_2-z_1\in L_+$. Consequently $f(z_1)\leqslant f(z_2)$ for all $z_1\in B^-$ and $z_2\in B^+$, i.~e. $p^-\leqslant p^+$. Take a number $p$, $p^-\leqslant p \leqslant p^+$ and put
\[
f^{\prime}(y+\lambda x_0)=p+\lambda f(x_0).
\]

In such a way we define an extension $f^{\prime}$ of $f$ from $B$ on $B^{\prime}$. From the definition directly implies that $f^{\prime}$ is a weakly additive functional. We will show that $f^{\prime}$ is order-preserving. It is order-preserving on $B$ owing to $f^{\prime}|B=f$. Besides it is evident that $f^{\prime}$ is order-preserving on $\{y+\lambda x_0\colon \ \lambda\in \mathbb{R}\}$. Let now $z-(y+\lambda x_0)\in L_+$, where $z\in B$. Then $(z-\lambda x_0)-y\in L_+$, i.~e. $(z-\lambda x_0)\in B^+$. That is why
\[
f^{\prime}(z-\lambda x_0) = f(z-\lambda x_0) = f(z)-\lambda f(x_0) \geqslant p^+ \geqslant p = f^{\prime}(y),
\]
i.~e. $f^{\prime}(z)\geqslant f^{\prime}(y+\lambda x_0)$. In the case when $(y+\lambda x_0)-z \in L_+$ one can similarly show that $f^{\prime}(z)\leqslant f^{\prime}(y+\lambda x_0)$.

Thus, a weakly additive, order-preserving continuous functional $f\colon B\to \mathbb{R}$ defining on an $A$-subspace $B$ can be extended to a weakly additive, order-preserving continuous functional $f^{\prime}\colon B^{\prime}\to \mathbb{R}$ on a wider $A$-subspace $B^{\prime}$ of $L$. At the same time the equality $f^{\prime}(x_0)=f(x_0)$ holds.

Consider the set of all pairs $(B^{\prime}, f^{\prime})$ such that $B\subset B^{\prime}\subset L$ where $B^{\prime}$ is an $A$-subspace, $f^{\prime}\colon B^{\prime}\to \mathbb{R}$ is a weakly additive, order-preserving continuous extension of $f$. The relation $(B^{\prime}, f^{\prime})\leqslant (B^{\prime\prime}, f^{\prime\prime})$ meaning that $f^{\prime\prime}\colon B^{\prime\prime}\to \mathbb{R}$ is a weakly additive, order-preserving continuous extension of $f^{\prime}$ on a subspace $B^{\prime\prime}$, $B^{\prime}\subset B^{\prime\prime}\subset L$, turns this set into a partially ordered set in which all chains are bounded. By Zorn's lemma there is the maximal element $(B_0, f_0)$ of this set. We will show that $B_0=L$.

Suppose that $B_0\neq L$. Take any point $y\in L\setminus B_0$ and put $B_1=B_0\cup \{y+\lambda x_0\colon \ \lambda\in \mathbb{R}\}$. Then $f_0$ can be extended to $f_1\colon  B_1\to \mathbb{R}$, and consequently, $(B_0, f_0)\leqslant (B_1, f_1)$. We got a contradiction with maximality of $B_0$. So, $B_0=L$.

\end{Proof}

\section{Uniform Boundedness Principle for Order-Preserving Operators}

Let $(E, \leqslant)$ and $(F, \leqslant)$ be partially ordered vector spaces.

\begin{Def}
A map $T\colon E\to F$ is said to be an order-preserving operator if for arbitrary points $x, y\in E$ the inequality $x\leqslant y$ implies $T(x)\leqslant T(y)$.
\end{Def}

Let $(E, \leqslant)$ be a partially ordered vector space with an order unit $1_E$ and $(F, \leqslant)$ be a partially ordered vector space.

\begin{Def}
A map $T\colon E\to F$ is said to be a weakly additive operator if $T(x+\lambda 1_E)=T(x)+\lambda T(1_E)$ takes place for each $x\in E$ and $\lambda\in \mathbb{R}$.
\end{Def}

The last definition immediately implies $T(0_E)=T(1_E-1_E)=T(1_E)-T(1_E)=0_F$, i.~e. $T(0_E)=0_F$ for a weakly additive operator $T\colon E\to F$.

The following statement shows weakly additive, order-preserving operators of vector spaces with an order unit are automatical continuous.

\begin{Prop}\label{opercont}
If $E$ and $F$ are partially ordered linear topological spaces with an order unit then each weakly additive, order-preserving operator $T\colon E\to F$ is continuous.
\end{Prop}

\begin{Proof}
 We will show the operator $T$ is continuous at zero $0_E$. At first we note the following case. If $T(1_E)=0_F$ then $T(x)=0_F$ for all $x\in E$ since $T$ is a weakly additive and order-preserving operator, and for every $x$ there exists $\lambda>0$ such that $-\lambda 1_E\leqslant x \leqslant \lambda 1_E$. So $T(E)\subset \{0_F\}$. This case we will not consider, i.~e. suppose $T(1_E)\neq 0_F$. Then $\|T(1_E)\|\neq 0$.

Let $V(0_F, \varepsilon)=\{y\in F\colon\ -\varepsilon 1_F<y<\varepsilon 1_F\}$
be a neighbourhood of zero $0_F$ in $F$, where $\varepsilon> 0$. Take the neighbourhood
$U\left(0_E, \frac{\varepsilon}{\|T(1_E)\|}\right)$ of zero $0_E$ in $E$.  For each vector $x\in U$ we have $-\frac{\varepsilon}{\|T(1_E)\|}1_E<x<\frac{\varepsilon}{\|T(1_E)\|}1_E$. Then $-\frac{\varepsilon}{\|T(1_E)\|}T(1_E)<T(x)<\frac{\varepsilon}{\|T(1_E)\|}T(1_E)$ since $T$ is a weakly additive, order-preserving operator. From here we get $\|T(x)\|<\varepsilon$, i.~e. $T(U)\subset V$. Thus $T$ is continuous at $0_E$. The following statement will finish the Proof.

\begin{Prop}\label{oper0cont}
If a weakly additive, order-preserving operator $T\colon E\to F$ of spaces with an order unit is continuous at zero then it is continuous on all $E$.
\end{Prop}

The proof of this Proposition is similarly to the Proof of Proposition \ref{continuous}.

\end{Proof}

\begin{Rem}\label{T1inner}
{\rm It is obvious that each linear non-negative operator on spaces with an order unit is weakly additive and order-preserving. The converse, in general, is not true. But, nevertheless, such operators are linear on a one-dimensional subspace $\{\lambda 1_E\colon\ \lambda\in \mathbb{R}\}\subset E$. In this case the image of the subspace $\{\lambda 1_E\colon\ \lambda\in \mathbb{R}\}$ at the map $T$ is, as clearly, a one-dimensional subspace $\{\lambda T(1_E)\colon\ \lambda\in \mathbb{R}\}\subset F$. We have $T(1_E)\in F_+$ but it is optional $T(1_E)\in Int F_+$. Therefore $T(1_E)$ is an order unit in $T(E)$ but it is optional to be an order unit in $F$. From here and  Proposition \ref{opercont} follows that for every weakly additive, order-preserving operator $T\colon E\to F$ on spaces $E$, $F$ with an order unit the inequality $\|T(1_E)\|<\infty$ takes place.}
\end{Rem}

Remind the following notions. A set $A$ in a normed space $E$ is called \textit{bounded} if there exists $R>0$ such that $A$ can be placed into the ball $\{x\in E\colon\ \|x\|\leqslant R\}$. A map $T\colon E\to F$ of normed spaces is called \textit{bounded} if it carries over a bounded set in $E$ to a bounded set in $F$. It is obvious that the boundedness of the map $T$ is equivalent to limitation of the set $\{\|T(x)\|\colon\ x\in E, \|x\|\leqslant R\}$ for every $R>0$. In other words, $\sup\{\|T(x)\|\colon\ x\in E, \|x\|\leqslant R\}<\infty$ for every bounded map $T$ and for each $R>0$.

The following statement shows weakly additive, order-preserving operators of vector spaces with an order unit are automatical bounded.

\begin{Prop}
Each weakly additive, order-preserving operator $T\colon E\to F$ of spaces with an order unit is bounded.
\end{Prop}

The proof follows from Remark \ref{T1inner}.

Let $E$ and $F$ be vector spaces with an order unit, $1_E$ and $1_F$, respectively. A collection $\mathcal{H}$ of weakly additive, order-preserving operators $T\colon E\to F$ is said to be \textit{equicontinuous} if to every neighbourhood $V$ of zero in $F$ there corresponds a neighbourhood $U$ of zero in $E$ such that $T(U)\subset V$ for all $T\in \mathcal{H}$. If the collection $\mathcal{H}$ consists only one weakly additive, order-preserving operator $T$, then $\mathcal{H}$ is equicontinuous as $T$ is continuous, and $\mathcal{H}$ is uniform bounded owing to boundedness of $T$. The following statement shows that each equicontinuous collection of weakly additive, order-preserving operators on vector spaces with an order unit is uniform bounded.

\begin{Prop}\label{unibound}
Let $E$ and $F$ be vector spaces with an order unit, $\mathcal{H}$ an equicontinuous collection of weakly additive, order-preserving operators $T\colon E\to F$, and $A$ a bounded subset of $E$. Then for every $T\in \mathcal{H}$ there exists a bounded subset $B$ of $F$ such that $T(A)\subset B$.
\end{Prop}

\begin{Proof}
Put $B=\bigcup\limits_{T\in \mathcal{H}}T(A)$. Since the collection $\mathcal{H}$ is equicontinuous then for every neighbourhood $V=V(0_F, \varepsilon)$ of zero in $F$ there exists a neighbourhood $U=U(0_E, \delta)$ of zero in $E$ that $T(U)\subset V$ for all $T\in\mathcal{H}$. So far as $A$ is bounded for enough big $t\in \mathbb{R}$ we have $A\subset tU$. It is clear, that $T(A)\subset T(t U)$. Assume that $x\in tU$. Then $\|x\|<t\delta$, i.~e. $-t\delta<x<t\delta$. As $T$ is weakly additive and order-preserving we have $-t\delta T(1_E)<T(x)<t\delta T(1_E)$, $\|T(x)\|<t\delta \|T(1_E)\|$, consequently,  $\|\frac{1}{t}T(x)\|<\delta \|T(1_E)\|=\|T(\delta 1_E)\|\leqslant \varepsilon$. Hence, $T(t U)\subset t V$. Thus $T(A)\subset t V$ for all $T\in \mathcal{H}$. It means that $B\subset t V$, i. e. the set $B$ is bounded.

\end{Proof}

The following result is a weakly additive, order-preserving operators' variant of the Banach-Steinhaus theorem.

\begin{Th}\label{Banach-Steinhaus}
Let $E$ and $F$ be vector spaces with an order unit, $\mathcal{H}$ be a collection of weakly additive, order-preserving operators $T\colon E\to F$, and $A$ be a set consisting of such points $x\in E$ that each orbit $\mathcal{H}(x)=\{T(x)\colon T\in \mathcal{H}\}$ is bounded in $F$. If $A$ is a set of the second category then $A=E$ and the collection $\mathcal{H}$ is equicontinuous.
\end{Th}

\begin{Proof}
Let $V=V(0_F, \varepsilon)$ and $W=W(0_F, \varepsilon^{\prime})$ be neighbourhoods such that $\overline{V}+\overline{V}\subset W$ where $\overline{V}$ is the closure of $V$ with respect to order topology in $F$. Put $B=\bigcap\limits_{T\in \mathcal{H}}T^{-1}(\overline{V})$. Let $x\in A$. Then for some positive integer $n$ we have $\mathcal{H}(x)\subset n V$ by virtue of boundedness of $\mathcal{H}(x)$. Hence $T(x)\in n V$ or $x\in n T^{-1}(V)$ for all $T\in \mathcal{H}$. It means that $x\in n B$. Thus $A\subset\bigcup\limits_{n=1}^{\infty} n B$. Thence at least one of the sets $n B$ is the second category owing to $A$ is so. A map $x\mapsto n x$ is a homeomorphism $E$ onto itself. Consequently the set $B$ is the second category. Continuity of operators $T\in \mathcal{H}$ implies $B$ is closed in $E$. As $B$ is the second category set, it has an inner point. By the construction of $B$ one can see that $\delta 1_E$ lies in $B$ as an inner point for enough small $\delta\in \mathbb{R}$. Let $\delta 1_E$ be such an inner point in $B$. Then a set $B-\delta 1_E=\{x-\delta 1_E\colon x\in B\}$ contains some neighborhood $U=U(0_E, \delta^{\prime})$ of zero and
\begin{multline*}
T(U)\subset T(B-\delta 1_E) = \{T(x-\delta 1_E)\colon x\in B\} =\\
=\{T(x)-\delta T(1_E)\colon x\in B\} = T(B)-\delta T(1_E)\subset \overline{V}-\overline{V}\subset W
\end{multline*}
for all $T\in \mathcal{H}$. It means that $\mathcal{H}$ is a equicontinuous collection. Then $\mathcal{H}$ is uniform bounded by Proposition \ref{unibound}. That is why an orbit $\mathcal{H}(x)$ is bounded for each $x\in E$.

\end{Proof}

If a vector space with an order unit is is a Banach space with respect to order norm then it said to be \textit{a complete space with an order unit}. As each Banach space is a set of the second category then Theorem \ref{Banach-Steinhaus} directly implies

\begin{Cor}\label{equicont}
Let $E$ be a complete space with an order unit and $F$ a vector space with an order unit, $\mathcal{H}$ a collection of weakly additive, order-preserving operators $T\colon E\to F$, and a collection $\mathcal{H}(x)=\{T(x)\colon T\in \mathcal{H}\}$ bounded in $F$. Then $\mathcal{H}$ is an equicontinuous collection.
\end{Cor}

As Proposition \ref{unibound} holds then Corollary \ref{equicont} means that a pointwise boundedness of an arbitrary collection weakly additive, order-preserving operators from a complete space with an order unit into a vector space with an order unit implies a uniform boundedness of this collection.

Let Let $E$ and $F$ be vector spaces with an order unit, $\{T_n\}$ a sequence of weakly additive, order-preserving operators $T_n\colon E\to F$. If for every $x\in E$ there exists a limit $\lim\limits_{n\to\infty}T_n(x)$ then putting
\begin{equation}\label{lim}
T(x)=\lim\limits_{n\to\infty}T_n(x),\qquad x\in E,
\end{equation}
we have a weakly additive, order-preserving operator. Indeed,
\[
T(x+\lambda 1_E)=\lim\limits_{n\to\infty}T_n(x+\lambda 1_E)=\lim\limits_{n\to\infty}(T_n(x)+\lambda T_n(1_E)) = T(x)+\lambda T(1_E),
\]
and if $x\leqslant y$ then
\[
T(x)=\lim\limits_{n\to\infty}T_n(x)\leqslant\lim\limits_{n\to\infty}T_n(y)=T(y).
\]

\begin{Cor}\label{limit}
Let Let $E$ and $F$ be vector spaces with an order unit, $\{T_n\}$ a sequence of weakly additive, order-preserving operators $T_n\colon E\to F$. If there exists a limit $\lim\limits_{n\to\infty}T_n(x)$, $x\in E$, then an operator $T\colon E\to F$ defined by (\ref{lim}) is also a weakly additive, order-preserving operator.
\end{Cor}

\section{Order-Preserving Variant of Open Mapping Theorem}

Remind that a map $f\colon X\to Y$ of topological spaces is called \textit{open} at $x_0\in X$ if for every open neighbourhood of $x_0$ in $X$ there exists an open neighbourhood $V$ of $f(x_0)$ in $Y$ such that $V\subset f(U)$. A map is \textit{open} on a topological space $X$ if it is open at every point of $X$.

\begin{Lem}\label{0open}
Let $E$ and $F$ be vector spaces with an order unit, $T\colon E\to F$ a weakly additive, order-preserving onto operator. If $T$ is open at zero then it is open on all $E$.
\end{Lem}

\begin{Proof}
Let for every neighbourhood $U=U(0_E, \varepsilon)$ of $0_E$ its image $T(U)=\{T(x): x\in U\}$ be open. We have $0_F\in T(U)$ as $T(0_E)=0_F$. Thence there exists an open neighbourhood $V=V(0_F, \delta)$ of $0_F$ such that $V\subset T(U)$.

Now let $x_0\in E$ be an arbitrary point and $U(x_0, \varepsilon$ a neighbourhood of $x_0$ got by shifting $U(0_E, \varepsilon)$ on vector $x_0$. Besides let $V(T(x_0), \delta)$ be a neighbourhood of $T(x_0)$ got by shifting $V(0_F, \delta)$ on vector $T(x_0)$. The proof of the Lemma will finished if we show that the following diagram is true

\[
\begin{matrix}
y\in V(T(x_0), \delta) & \overset{(1)}\Longleftrightarrow & y-T(x_0)\in V(0_F, \delta)\\
&  & \mbox{\scriptsize (2)}\Downarrow\\
y\in T(U(x_0, \varepsilon)) & \overset{(3)}\Longleftrightarrow & y-T(x_0)\in T(U(x_0, \varepsilon)).
\end{matrix}
\]

The equivalence of the double inequalities $-\delta 1_F < y-T(x_0) < \delta 1_F$ and $T(x_0) - \delta 1_F < y < \delta 1_F + T(x_0)$ implies $(1)$. Since $V\subset T(U)$ we have $(2)$. And the equivalence of the double inequalities $-\varepsilon T(1_E) < y-T(x_0) < \varepsilon T(1_E)$ and $T(x_0) - \varepsilon T(1_E) < y < \varepsilon T(1_E) + T(x_0)$ implies $(3)$.

Thus for an arbitrary point $x\in E$ and its arbitrary neighbourhood $U=U(x, \varepsilon)$ there exists open neighbourhood $V=V(T(x), \delta)$ such that $V\subset T(U)$.

\end{Proof}

Since order topology is invariant with respect to the shift of points of the vector space Lemma \ref{0open} implies

\begin{Cor}
Weakly additive, order-preserving surjective operator of vector spaces with an order unit is open iff it is open at zero.
\end{Cor}

Remind that a metric $d$ on a vector space $E$ is invariant concerning to a shift of points of $E$ if $d(x+z, y+z)=d(x, y)$ for all $x, y, z\in E$. Define an order metric by the rule
\[
d(x, y)=\|y-x\|=\inf\{\lambda>0:\ -\lambda 1_E <y-x< \lambda 1_E\}.
\]

It is easy to see that the following assertion holds.

\begin{Lem}
The order metric on a vector space with an order unit is invariant according to a shift of points.
\end{Lem}

Let $E$ and $F$ be vector spaces with an order unit. A product $E\times F$ over $(0_E, 0_F)$ becomes a vector space with an order unit if we will introduce to it coordinatewise operations of sum and multiplication by number
\[
\alpha(x_1, x_2) + \beta(y_1, y_2) = (\alpha x_1+\beta y_1,\ \alpha x_2+\beta y_2),
\]
and coordinatewise partially order
\[
(x_1, x_2) \leqslant (y_1, y_2) \Leftrightarrow (x_1 \leqslant y_1 \mbox{ and } x_2 \leqslant y_2).
\]
Order norm on $E\times F$ is defined by the rule
\[
\|(x_1, x_2)\|=\inf\{\lambda>0: -\lambda(1_E, 1_F) \leqslant (x_1, x_2)\leqslant \lambda(1_E, 1_F)\}.
\]

Here $(1_E, 1_F)$ is one of inner points of $(E\times F)_+=E_+\times F_+$ that is why without losing generality we assume $(1_E, 1_F)$ is an order unit in the product. Denote $1_{E\times F}=(1_E, 1_F)$.

Let $T\colon E\to F$ be a weakly additive, order-preserving operator. The set of all pairs $(x, T(x))$, $x\in E$, is called \textit{a graph} of $T$.

\begin{Lem}\label{graph}
Let $E$ and $F$ be vector spaces with an order unit, $1_E$ an order unit in $E$, $T\colon E\to F$ a weakly additive, order-preserving operator. Then the graph $G$ of $T$ is an $A$-subspace of $E\times T(E)$ with an order unit $1_{E\times T(E)}$.
\end{Lem}

\begin{Proof}
We have $(0_E, 0_F)\in G\subset E\times T(E)$ since $T(0_E)=0_F$. Consider $(x_1, x_2)\in E\times T(E)$ and $\lambda \in \mathbb{R}$. Then
\begin{multline*}
(x_1, x_2) + \lambda 1_{E\times T(E)}=(x_1, T(x_1)) + (\lambda 1_E, \lambda 1_{T(E)}) =\\
=(x_1+\lambda 1_E, T(x_1)+\lambda 1_{T(E)}) = (x_1+\lambda 1_E, T(x_1) + \lambda T(1_E))=\\
= (x_1+\lambda 1_E, T(x_1 + \lambda 1_E)),
\end{multline*}
i.~e. $(x_1, x_2) + \lambda 1_{E\times T(E)}\in G$.

\end{Proof}

Lemma \ref{graph} and Remark \ref{T1inner} imply

\begin{Cor}
Let $E$, $F$ be vector spaces with an order unit, $1_E$ an order unit in $E$, $T\colon E\to F$ a weakly additive, order-preserving operator. Then the image $T(E)$ is $A$-subspace of $F$ if and only if $T(1_E)\in Int F_+$.
\end{Cor}

\begin{Rem}\label{1f=t1e}
{\rm Further, during current section, without losing of generality, we will consider such weakly additive, order-preserving operators $T$ for which $T(1_E)\in Int F_+$. Then we may assume that $T(1_E)$ is an order unit in $F$. Put $1_F=T(1_E)$.}
\end{Rem}

At last we will form a variant of the Open Mapping Theorem for weakly additive, order-preserving operators.

\begin{Th}\label{openmapth}
Let $E$ be a complete space with an order unit, $F$ a vector space with an order unit, and $T\colon E\to F$ a weakly additive, order-preserving operator such that $T(E)=F$ and $F$ is a set of the second category. Then

(i) the map $T$ is open;

(ii) $F$ is a complete space with an order unit.
\end{Th}

\begin{Proof}
Let $U(0_E, \varepsilon)$ be an open neighbourhood. Then according to Remark \ref{1f=t1e} we have
\begin{multline*}
T(U(0_E, \varepsilon))= \{T(x)\in F: -\varepsilon 1_E <x< \varepsilon 1_E\}=\\
=\{T(x)\in F: -\varepsilon 1_F <T(x)< \varepsilon 1_F\}=(\mbox{by condition})=\\
=\{y\in F: \mbox{ there exists } x\in U(0_E, \varepsilon) \mbox{ such that } y=T(x) \mbox{ and }
-\varepsilon 1_F <y< \varepsilon 1_F\}=U(0_F, \varepsilon).
\end{multline*}

It reminds to show that $(ii)$ takes place.

Let $\{y_n\}\subset F$ be a fundamental sequence. Then for every $\varepsilon>0$ there exists a number $n$ such that at all $k\geqslant n$, $m\geqslant n$ the double inequalities
\[
-\varepsilon 1_F <y_m-y_k<\varepsilon 1_F
\]
hold. One may assume $\varepsilon=\frac{1}{n}$. Then $-\frac1n 1_F <y_m-y_k<\frac1n 1_F$. Since $T(U(0_E, \frac1n))=U(0_F, \frac1n)$ there exists $x_m, x_k\in E$ such that $T(x_m)=y_m$, $T(x_k)=y_k$ and $-\frac1n 1_E<x_m-x_k< \frac1n 1_E$. So we have constructed a fundamental  sequence $\{x_n\}\subset E$. By completeness of $E$ the sequence have a limit $x=\lim\limits_{n\to\infty}x_n$. As $T$ is continuous we have $T(x)=\lim\limits_{n\to\infty}T(x_n)=\lim\limits_{n\to\infty}y_n$. Then $\lim\limits_{n\to\infty}y_n\in T(E)=F$. Thus, $F$ is complete space with an order unite.

\end{Proof}

\begin{Rem}
{\rm Note that open mapping principle for weakly additive, order-preserving operators it is impossible to form as the linear case. In the distinguishing from the linear case, weakly additivity and order-preserving of $T$, and being of $T(E)$ the second category set does not imply the equality $T(E)=F$. On the other hand $T(E)$ must not be open in $F$. At last if $T$ is not onto in Lemma \ref{0open} then openness of $T$ at zero does not provide it openness on all the space.

Note that in linear topological spaces does not exist open subspace without all space. But an $A$-subspace, distinguished from the subspace, may be open, closed or everywhere dense in the vector space with order unit.}
\end{Rem}

\begin{Ex}
{\rm Let $L=\mathbb{R}^2=\{(x_1, x_2): \ x_i\in \mathbb{R},\ i=1,\ 2\}$ be vector space with an order unit considered in Example \ref{positbutnotord}. Then  $L_+=\{(x_1, x_2)\in \mathbb{R}^2: \ x_i\geqslant 0, i=1,\ 2\}$ is a positive cone in $L$. Fix $\mathbf{1}=(1, 1)\in Int\mathbb{R}^2_+=\{(x_1, x_2)\in \mathbb{R}^2: x_1>0, x_2>0\}$ as an order unit in it.

a) It is easy to see that the set $B=\{(x_1, x_1+a)\in \mathbb{R}^2: -1<a<1\}$ is an open (with respect to order topology) $A$-subspace, but $B\neq \mathbb{R}^2$.

b) Let $\mathbb{Q}$ be the set of rational numbers. Then $C=\{(x_1, x_1+r)\in \mathbb{R}^2; r\in \mathbb{Q}\}$ is an everywhere dense $A$-subspace in $\mathbb{R}^2$.

c) The set $D=\{(x_1, x_1+a)\in \mathbb{R}^2: -1\leqslant a \leqslant 1\}$ is a closed $A$-subspace in $\mathbb{R}^2$.

d) Define a map $T\colon \mathbb{R}^2\to \mathbb{R}^2$ by the rule
\[
T(x_1, x_2)=\begin{cases}
(x_1, x_1-1), & \mbox{ at } \quad x_2\leqslant x_1-1,\\
(x_1, x_2), & \mbox{ at } \quad x_1-1< x_2 < x_1+1, \\
(x_1, x_1+1), & \mbox{ at } \quad x_2\geqslant x_1+1.
\end{cases}
\]
It is easy to check that $T$ is a weakly additive map. Let us show that the map $T$ is order-preserving. It clear that $T$ is order-preserving on $B$ by $T=id_B$.

Let $x_2\geqslant x_1+1$. Take a vector $(y_1, y_2)\in \mathbb{R}^2$ such that $(x_1, x_2)\leqslant (y_1, y_2)$. The following three cases possible.

Case 1) $y_2\geqslant y_1+1$. Then
\[
T(x_1, x_2) = (x_1, x_1+1)\leqslant (y_1, y_1+1) = T(y_1, y_2).
\]

Case 2) $y_1-1\leqslant y_2 \leqslant y_1+1$. Then $x_1+1\leqslant y_2$. That is why
\[
T(x_1, x_2) = (x_1, x_1+1)\leqslant (y_1, y_2) = T(y_1, y_2).
\]

Case 3) $y_2\leqslant y_1-1$. Then $x_1+1\leqslant y_1-1$. Censequently
\[
T(x_1, x_2) = (x_1, x_1+1)\leqslant (y_1, y_1-1) = T(y_1, y_2).
\]

Similarly, one may show that $T$ is order-preserving when $x_2\leqslant x_1-1$. Thus $T$ is order-preserving on all $\mathbb{R}^2$.

We have $T(\mathbb{R}^2)=D\neq \mathbb{R}^2$ though the operator $T$ is weakly additive and order-preserving, and the image $T(\mathbb{R}^2)$ is the second category. Clearly the image $T(\mathbb{R}^2)$ is closed in $\mathbb{R}^2$ and it is not open. Moreover $T$ is open at zero but it is not open on $\mathbb{R}^2$. Really for the open neighbourhood $U((2, 4), 1)=\{(x_1, x_2)\in \mathbb{R}^2: 1<x_1<3, 3<x_2<5\}$ of the point $(2, 4)\in \mathbb{R}^2$ its image $T(U)=\{(x_1, x_1+1): 1<x_1<3\}$ is not open in $T(\mathbb{R}^2)$.}
\end{Ex}

\section{Order-Preserving Variant of Banach-Alaoglu Theorem}

Let $E$ be a vector space with an order unit. Fix $1_E$ as an order unit. By $E_+^W$ we the set of all weakly additive, order-preserving functionals $f\colon E\to \mathbb{R}$. On $E_+^W$ define algebraic operations pointwise. Then $E_+^W-E_+^W$ turns to a vector space with an order unit. Denote $E^W=E_+^W-E_+^W$. Put $E^O=\{f\in E_+^W: f(1_E)=1\}$. Provide $E^W$ with the pointwise convergence topology. A collection of the sets of the view
\[
\langle f; x_1,\dots, x_n; \varepsilon\rangle=\{g\in E^W: |f(x_i)-g(x_i)|<\varepsilon,\ i=1,\dots, n\}
\]
forms a base of open neighbourhoods of $f\in E^W$, where $\varepsilon>0$, $x_i\in E$, $i=1,\dots, n$.

The main result of the section is the following variant of the Banach-Alaoglu theorem for weakly additive, order-preserving functionals.

\begin{Th}\label{Banach-Alaoglu}
If $V$ is a neighbourhood of zero in $E$ then the set
\[
K=\{f\in E^O: |f(x)|\leqslant 1\}
\]
is a compact in the pointwise convergence topology.
\end{Th}

\begin{Proof}
Since neighbourhoods of zero are absorbing sets, for every point $x\in E$ there exists $\gamma(x)\ \mathbb{R}_+$ such that $x\in \gamma(x)V$. That is why $|f(x)|\leqslant \gamma(x)$ for all $f\in E^W$ and $x\in E$. For every $x\in E$ denote $D_x=[-\gamma(x), \gamma(x)]$ and assume that $\tau$ is the Tychonoff topology in the product $P=\prod\limits_{x\in E}D_x$. It is well known that $P$ is a Hausdorff compact space. By the construction we have $K\subset P\cap E^W$. We will show that $K$ is closed in $P$. Let $f_0\in P$ and $f_0=f_0^+-f_0^-$, where $f_0^+, f_0^-\in P\cap E_+^W$. Suppose $\{f_{\alpha}^+\}\subset P\cap E_+^W$ and  $\{f_{\theta}^+\}\subset P\cap E_+^W$ are nets converging to $f_0^+$ and $f_0^-$ respectively. Then owing to Corollary \ref{limit} we have $f_0^+, f_0^-\in E^W$, and therefore $f_0\in E^W$. On the other hand $|f_0(x)|=|(f_0^+(x)-f_0^-(x))|\leqslant \max\{|f_0^+(x)|, |f_0^-(x)|\}\leqslant \gamma(x)$ by $|f_{\alpha}^+(x)|\leqslant \gamma(x)$ and $|f_{\alpha}^-(x)|\leqslant \gamma(x)$ for all $x\in E$, $\alpha$ and $\theta$. Therefore $|f_0^+(x)|\leqslant \gamma(x)$ for all $x\in E$ and $|f_0^+(x)|\leqslant 1$ so far as $x\in V$. It means that $f_0\in K$.

\end{Proof}

\begin{Cor}
$E^O$ is a compact in the pointwise convergence topology.
\end{Cor}

If $E$ is a separable vector space with an order unit then Theorem \ref{Banach-Alaoglu} improves as

\begin{Th}
If $E$ is a separable vector space with an order unit, and $K$ is a compact (with respect to pointwise convergence topology) subspace of $E^W$ then $K$ is metrizable.
\end{Th}

\begin{Proof}
Let $\{x_n\}$ be countable everywhere dense subset of $E$. For every $f\in E^W$ put $M_n(f)=f(x_n)$. By the definition of pointwise convergence topology every $M_n$ is a continuous function on $E^W$. If $M_n(f)=M_n(f^\prime)$ for all $n$ then continuous functions $f$ and $f^\prime$ coinside on everywhere dense subset. Thus $\{M_n\}$ is a countable family of continuous functions which separate points of the space $E^W$, in particular of $K$. Hence $K$ is metrizable as each Hausdorff compact space which has a countable sequence of real-valued functions separating its points is metrizable.

\end{Proof}

\begin{Cor}
If $E$ is separable vector space with an order unit then $E^O$ is a metrizable compact in the pointwise convergence topology.
\end{Cor}


\begin{thebibliography}{00}


\bibitem{grabroub2000appChoquetint}
Grabisch~M. and Roubens~M. \textsl{Application of the Choquet Integral in Multicriteria Decision Making},  In M. Grabisch, T. Murofushi, and M. Sugeno, editors, Fuzzy Measures and Integrals -- Theory and Applications, pages 348–374. Physica Verlag, 2000.

\bibitem{radul2018idemmeasurAR}
Radul~T. \textsl{Idempotent Measures: Absolute Retracts And Soft Maps}, arXiv:1810.09140v1 [math.GN] 22 Oct 2018

\bibitem{vern2009vectspordunit}
Vern~I. Paulsen and Mark Tomforde, \textsl{Vector Spaces with an Order Unit},
arXiv:0712.2613v4 [math.OA] 10 Jun 2009

\bibitem{pep2018unibound}
Aljo\v{s}a Peperko, \textsl{Uniform Boundedness Principle for Nonlinear Operators on Cones of Functions}, Hindawi Journal of Function Spaces, Volume 2018, Article ID 6783748

\bibitem{zai2006openmap}
Zaitov~A.~A.\textsl{On extension of order-preserving functionals}, Reports of Academy of Sciences of the Republic of Uzbekistan. 2005, No 5. P. 3-7.

\bibitem{zai2005BanachAlaoglu}
Zaitov~A.~A. \textsl{Banach-Alaoglu theorem for order-preserving functionals},
Theses of reports of the international scientific conference "Operator Algebras and Quantum Probability Theory", September 7-10, 2005 P. 81-83.

\bibitem{zai2006openmap}
Zaitov~A.~A. \textsl{Open mapping theorem for order-preserving operators}, The collection of theses of the International conference of young scientists devoted 1000 to the anniversary of Mamun Academy of Khwarezm. Tashkent, 2006, P. 4







\end{thebibliography}
\end{document}